\documentclass{amsart}

\usepackage{graphicx}
\graphicspath{ {images/} }
\usepackage{enumitem} 
\usepackage{amssymb} 
\usepackage{bbold} 

\graphicspath{ {images/} }

\newtheorem{theorem}{Theorem}[section]
\newtheorem{lemma}[theorem]{Lemma}

\theoremstyle{definition}
\newtheorem{definition}[theorem]{Definition}

\newtheorem{conj}{Conjecture}
\newtheorem{prop}{Proposition}

\theoremstyle{remark}
\newtheorem{remark}[theorem]{Remark}

\numberwithin{equation}{section}

\title{MaxDim of some Simple groups}
\author{Liu Tianyue}
\author{R. Keith Dennis}
\date{July 2017}

\begin{document}

\maketitle

\section{Introduction}
There are various generalizations of the concept of dimension for vector spaces in finite group theory. In this paper we are primarily concerned with three of them, namely $m(G)$, $MaxDim(G)$ and $i(G)$.
\begin{definition}
A list of subgroups $\{H_1,...,H_n\}$ are said to be in general position if for all $i \in \{1,...,n\}$, $\displaystyle\bigcap\limits_{j\neq i} H_j\gneq \displaystyle\bigcap\limits_j H_j$.\par
\begin{enumerate}[label=(\roman*)]
\item$MaxDim(G)$ of a finite group $G$ is defined as the maximal size of a list of $G$'s maximal subgroups in general position.\par
\item$m(G)$ is defined as the maximal size of an irredundant generating sequence in $G$\par
\item$i(G)$ is defined as the maximal size of an irredundant sequence in $G$\par
\end{enumerate}
\end{definition}
The basic relation between those group invariants is captured by the following theorem.
\begin{theorem}
For a finite group $G$, $m(G)\leq MaxDim(G)\leq i(G)$
\end{theorem}
$MaxDim$, $m$ and $i$ can also be defined using the subgroup lattice. For a set $[n]$ with $n$ elements, the power set $P([n])$ has a partial order induced by inclusion. We shall call a poset $T$ a tower with dimension $n$ if it is isomorphic to $P([n])$, and call the pre-images of $(n-i)$-element sets in $P[n]$ the $i$-th level of the tower. In this context, we can see that the following proposition holds.
\begin{prop}
\begin{enumerate}[label=(\roman*)]
\item $i(G)$ is the dimension of the largest tower embeddable in the subgroup lattice of $G$.
\item $MaxDim(G)$ is the dimension of the largest tower embeddable in the subgroup lattice of $G$, so that all first level elements of the tower are maximal subgroups. That is, for a first level element $S$, $A\geq S\Rightarrow A=S$ or $A=G$.
\item $m(G)$ is the dimension of the largest tower embeddable in the subgroup lattice of $G$, so that the join of all $(n-1)$-th level elements of the tower is $G$, and all $(n-1)$-th level elements are the joins of the $n$-th level element and cyclic.
\end{enumerate}
\end{prop}
Before we proceed to the results on the $MaxDim$ of simple groups, we shall state some useful lemma here.
\begin{lemma}
$m(G)\leq \max \{i(M)|M$ is maximal subgroup of $G\}$+1\par
$i(G)=\max\{m(G), \max \{i(M)|M$ is maximal subgroup of $G\}\}$
\end{lemma}

\section{The Suzuki Groups}
In this section we compute m and MaxDim for some Suzuki groups of Lie type. The main result was almost implied by a proof in the thesis of Whiston, although he did not state the claim in his paper.
\begin{prop}
If $q$ is a prime power of 2, then $Suz(q)$ has $MaxDim=m=3$.
\end{prop}
\begin{proof}
To prove the proposition, we need to use the following theorem of Suzuki \cite{Wh01}:
\begin{theorem}
$Suz(q)$ acts 2-transitively on $q^2+1$ points. The maximal subgroups of $Suz(q)$ fall in one of the conjugacy classes listed below:\par
\begin{enumerate}[label=(\roman*)]
\item Stabilizers of a point, $E_q.E_q:C_{q-1}$
\item Stabilizers of a two point subset, $D_{2(q-1)}$
\item $C_{q+\sqrt{2q}+1}:C_4$
\item $C_{q-\sqrt{2q}+1}:C_4$
\item $Suz(q_0)$ when $q=q_0^r$, $r$ is a prime number and $q_0$ is not 2
\end{enumerate}
\end{theorem}

The general strategy behind the proof of Prop. 1 is to show on one hand, since Suzuki groups do not satisfy the $B$-properties, $m\geq 3$; on the other hand, the structure of the maximal subgroups does not allow $MaxDim$ to exceed 3. We shall prove the upper bound for MaxDim by inspecting the intersections of those maximal subgroups.\par

\begin{prop}
$q$ is an odd power of 2, $A$, $B$ are maximal subgroups of $Suz(q)$. The intersection $I:=A\cap B$ is:
\begin{enumerate}
\item Both $A$, $B$ are of type $(iii)$ or $(iv)$, then $I$ is a subgroup of $C_4$
\item $A$ is of type $(ii)$ and $B$ of type $(iii)$ or $(iv)$, then $I$ is a subgroup of $C_2$
\item $A$ is of type $(i)$ and $B$ of type $(iii)$ or $(iv)$, then $I$ is a subgroup of $C_4$
\item Both $A$, $B$ are of type $(ii)$, then $I$ is a subgroup of $C_2$
\item Both $A$, $B$ are of type $(i)$, then $I$ is isomorphic to $C_{q-1}$
\item $A$ is of type $(i)$ and $B$ of type $(ii)$, then $I$ is isomorphic to either $C_{q-1}$ or a subgroup of $C_2$
\end{enumerate}
\end{prop}
To prove claim $(1) \sim(4)$ we need the following lemma:
\begin{prop}
If $F:= C_n\rtimes B$,  $G:= C_m\rtimes B'$ are distinct maximal in simple group $S$, $gcd(|B|,m)=gcd(|B'|,n)=1$, then the intersection $I=F\cap G$ is isomorphic to a common subgroup of $H$ and $H'$, where $H\cong B$ and $H'\cong B'$.
\end{prop}
\begin{proof}
Assume there is a nontrivial $x\in C_n$ lying inside the intersection $I$, then $x\in C_m$ (and vice versa) because any other element in $G$ would have order dividing $|B'|$, while as an element of $C_n$, the order of $x$ cannot divide anything relatively prime to $n$.\par
On one hand, if we assume $C_n$ is not contained in the intersection, then there must be a nontrivial element $z\in C_n\backslash I$. Since $G$ is maximal in $S$, $S=<G,\{z\}>$. But $z$ and $x$ lie in the same cyclic group $C_n$, so $z<x>z^{-1}=<x>$. Moreover, by the cyclicity of $C_m$, $<x>$ is normal in $G$. That is to say, $g<x>g^{-1}=<x>$ for any finite word $g$ using $z$ and elements in $G$. Thus $<x>$ must be a proper normal subgroup of $S$, which contradicts $S$ being a simple group. If $C_m$ is not contained in the intersection, we yield the same contradiction.\par
On the other hand, if $C_n,C_m\subseteq I$ then the normalizer of $K:=C_n\cap C_m$ must contain both $F$ and $G$, so $S=N_S(K)$, also contradicting the simplicity of $S$. Thus the intersection $I$ must intersect trivially with $C_n$ or $C_m$.\par
Write elements in $F$ as ordered pairs. Since $I$ intersects $C_n$ trivially, for any $y\in B$ there is at most 1 element $(x,y)\in I$, so the homomorphism $(x,y)\mapsto y$ is injective. Similarly, the homomorphism $(x',y')\in G\mapsto y'\in B'$ is also injective. $I$ is isomorphic to a common subgroup of $B$ and $B'$. 
\end{proof}
Thus, since both $q-1$ and $q\pm\sqrt{2q}+1$ are odd, we have that $gcd(q\pm\sqrt{2q}+1,4)=1$ implies claim (1), $gcd(q-1,4)=gcd(q\pm\sqrt{2q}+1,2)=1$ implies claim (2), and $gcd(q-1,2)=1$ implies claim (4). Claim (3) follows from Whiston's result, namely since $gcd(4(q\pm\sqrt{2q}+1),q^2(q-1))=4$, the order of the intersection must divide 4 \cite{Wh01}. However, a Klein four group cannot be embedded into a subgroup of type $(iii)$ or $(iv)$, so the intersection must be a subgroup of $C_4$.\par
The intersection of two distinct point stabilizers, say $Stab(a)$ and $Stab(b)$, is equivalent to the stabilizer of $a$ in the pair-stabilizer $Stab(\{a,b\})$. The orbit of $a$ in $Stab(\{a,b\})$ is $\{a,b\}$, so by the Orbit-Stabilizer Theorem, the intersection of $Stab(a)$ and $Stab(b)$ has index 2 in dihedral group $Stab(\{a,b\})\cong D_{2(q-1)}$. Thus the intersection is cyclic of order $q-1$, and claim (5) is proven.\par
If $I$ is the intersection of $Stab(a)$ and $Stab(\{u,v\})$, then there are two cases. If $a\in\{u,v\}$ then $I\cong Stab(u)\cap Stab(v)\cong C_{q-1}$. If $a\notin\{u,v\}$, then we need the following lemma of Huppert and Blackburn \cite{HuBl82}:
\begin{prop}
Any nontrivial element of $Suz(q)$ has at most 2 fixed points.
\end{prop}
That is to say, the intersection of any three point stabilizer must be trivial. If $a\notin\{u,v\}$, then since $Stab(u)\cap Stab(v)\cap Stab(a)$ has index at most 2 in $I:=Stab(a)\cap Stab(\{u,v\})$, $I$ must be a subgroup of $C_2$ and claim (6) holds.\par
Since any subgroup of $C_4$ or $C_2$ has $i\leq1$, a length 4 list $l$of maximal subgroups in general position must include no subgroup of type $(iii)$ or $(iv)$, at most 1 subgroup of type $(ii)$, and at most 3 point stabilizers. Since 3-fold intersection of point stabilizers is trivial, they cannot coexist with another dihedral group in the list, so there cannot be a such a length 4 list $l$. That is to say,  $3=m=MaxDim=3$.
\end{proof}

\section{The Suzuki Tower}
The Suzuki tower is a family of five simple groups: $PSL_2(7)<PSU_3(3)<J_2<G_2(4)<Suz$. Here $J_2$  is the Hall-Janko group, $G_2(4)$ is an exceptional group of Lie type, and $Suz$ is the sporadic Suzuki group. Each of them act transitively on a distance transitive graph, denoted by $\Gamma_1$ to $\Gamma_5$ respectively. Moreover, for $i=1,2,3,4$, $\Gamma_i$ are the local subgraphs of $\Gamma_{i+1}$. That is to say, for any vertex $v\in \Gamma_{i+1}$, the subgraph induced by $v$'s neighbors is $\Gamma_i$. When the groups in the Suzuki tower acts as automorphisms of these graphs, the vertex stabilizers are isomorphic to the group one level beneath it in the tower \cite{Wi09}. Using these information and some knowledge on the structure of $\Gamma_1$, we can give an inductive argument to bound the $MaxDim$ of these groups from below.\par
\begin{prop}
The $MaxDim$ of $PSU_3(3)$, $J_2$, $G_2(4)$, and $Suz$ are greater or equal to 4, 5, 6 and 7 respectively. 
\end{prop}
\begin{proof}
As mentioned by Brouwer, $\Gamma_1$ is the co-Heawood graph on 14 vertices, which is the point-line non incidence graph of the Fano plane. In other words, it is the distance-3 graph of the Heawood graph in Figure 1 \cite{Br89}. $PSL_2(7)$ acts on it with rank 4, and the stabilizer of a vertex $v_3$ (e.g. the one labeled with a square in Figure 1) is $S_4$, which has vertex orbit size $1+4+6+3$. The 6 orbit corresponds to the vertices at distance 2 with $v_3$. Stabilizing another point in $v_3$'s neighborhood (those vertices labeled with circles in Figure 1) gives $S_3$, and stabilizing a point in the 6 orbit yields the Klein 4 group $C_2\times C_2$. Stabilizing 2 points in $v_3$'s neighborhood would result in $C_2$, because the one vertex stabilizer $S4$ acts on the four neighbors as the permutation group $S_4$. Notice that $v_3$'s neighbors are at distance 2 with each other in the co-Heawood graph.\par
Now to bound MaxDim of a Suzuki tower group at level $i$, $i>1$, we can pick an arbitrary vertex $v_{i+2}$ first. Then we keep picking vertices in the common neighborhood of previously picked vertices, until the common neighborhood becomes the co-Heawood graph. Since each time we add a new point the common neighborhood induces a graph one level beneath the previous graph in the tower, we should have collected $i-1$ vertices so far. Pick a vertex in the co-Heawood graph and denote it $v_3$. then pick two more vertices in $v_3$'s neighborhood and denote them by $v_2$, $v_1$.\par
I claim that the point stabilizers of $\{v_1,...,v_{i+2}\}$ in $\Gamma_i$ are maximal subgroups in general position. The intersection of all those point stabilizers is the same as the point-wise stabilizer of $v_1$, $v_2$, and $v_3$ in the co-Heawood graph, so it must be $C_2$.  Meanwhile, the subgraph of $\Gamma_i$ induced by $\{v_1,...,v_{i+2}\}$ is a complete $i+2$ graph with one edge removed, and its subgraphs induced on any $i+1$ vertices are complete $i+1$ graphs with or without an edge removed. If the subgraph on $i+1$ vertices is a complete graph, then the intersection of corresponding stabilizers is isomorphic to $S_3$, i.e. the stabilizer of a vertex and one of its neighbor in the co-Heawood graph. Otherwise, the instersection of point stabilizers must be $C_2\times C_2$, which is the point-wise stabilizer of two vertices at distance 2 in the co-Heawood graph. In either case, the $i+1$-fold intersections contains the overall intersection properly. thus by definition the $i+2$ point stabilizer are in general position.
\begin{figure}
\includegraphics[width=0.5\textwidth]{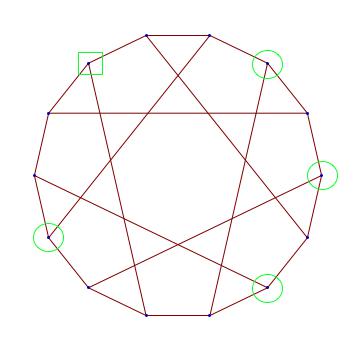}
\caption{Heawood Graph}
\end{figure}
\end{proof}
\begin{remark}
We actually can prove a stronger result here, namely the lower bounds for $MaxDim$ are the lower bonds for$m$ as well. Let $\{v_1,...,v_{i+2}\}$ be picked as mentioned above, and denote the two neighbors of $v_3$ not in the list by $a$ and $b$. Then there are two elements $\sigma,\tau$ in the point-wise stabilizer $Stab(\{v_3,...,v_{i+2}\})\cong S_4$ so that $\sigma: a\mapsto b,b\mapsto v_1,v_1\mapsto a, v_2\mapsto v_2$ and $\tau: a\mapsto v_2,v_2\mapsto a, b\mapsto b,v_1\mapsto v_1$. $\sigma$ and $\tau$ generates an $S_4$, which is exactly the point-wise stabilizer $Stab(\{v_3,...,v_{i+2}\})$. Thus if we pick any $\alpha_3,...,\alpha_{i+2}$ so that $\alpha_j\in \displaystyle\bigcap\limits_{k\neq j} Stab(v_k)\backslash \displaystyle\bigcap\limits_k Stab(v_k)$, then since for all $j\geq 2$, $Stab(\{v_j,...,v_{i+2}\})$ is maximal in $Stab(\{v_{j+1},...,v_{i+2}\})$, we know that $l:=\{\sigma,\tau,\alpha_3,...,\alpha_{i+2}\}$ generates the entire group. Since sequence $l$ certifies a list of maximal subgroups in general position, it must be irredundant, so we know $m\geq i+2$ for the group at level $i$. In particular, since the intersection of the point stabilizers is $C_2$, applying a lemma of Keith Dennis tells us $PSU_3(3)$ and $J_2$ fails the replacement property.
\end{remark}

\begin{remark}
The lower bounds of $MaxDim$ are strict for $PSU_3(3)$ and $J_2$, which has been checked computationally.
\end{remark}

\section{The Mathieu Groups}
In this section we give a computer-aided proof for the $MaxDim$ of Mathieu groups.
\begin{prop}
The $MaxDim$ of $M_{11}$, $M_{12}$, $M_{22}$, $M_{23}$, $M_{24}$ are 5, 6, 6, 6, and 7 respectively.
\end{prop}
\begin{proof}
$MaxDim$ for $M_{11}$ and $M_{12}$ has already been calculated by Brooks \cite{Br13}. \par
$M_{22},M_{23},M_{24}$ can be defined as automorphism groups of Steiner systems. A Steiner system $(a,b,v)$ is a set $\mathcal{S}$ on $v$ points with a collection of subsets $\mathcal{B}$ called "blocks", so that any block in $\mathcal{B}$ contains exactly $b$ points and any $a$ points are contained in exactly one block \cite{Wi09}. In this section We shall utilize the axioms of Steiner systems to prove our claim on $MaxDim$.\par
We know that $M_{24}$ is the automorphism group of the unique Steiner (5,8,24) system, and is 5-transitive on the 24 points. $M_{23}$ is the stabilizer of a point in $M_{24}$. As will be shown later, we can find a list of 7 point stabilizers in $M_{24}$ that are in general position. If $\{Stab(a_1),...,Stab(a_7)\}$ are in general position in $M_{24}$, then $\{Stab(a_1)\cap Stab(a_2),...,Stab(a_1)\cap  Stab(a_7)\}$ are in general position in $Stab(a_1)\cong M_{23}$. Since $M_{23}$ is 4-transitive on 23 points, point stabilizers are maximal in $M_{23}$. Thus $\{Stab(a_1)\cap Stab(a_2),...,Stab(a_1)\cap  Stab(a_7)\}$ is a list of maximal subgroups in general position, so $MaxDim(M_{23})\geq 6$. Applying the same argument to $M_{22}$ the point stabilizer in $M_{23}$, we know that $MaxDim(M_{22})\geq 5$.\par
However, as suggested by Sophie Le's computational results, for $M_{22}$ there is a way to construct a longer list of maximal subgroups in general position. Pick 2 points $\{a,b\}$ in Steiner $(5,8,24)$ system, the 22 points leftover together with the blocks that contain $\{a,b\}$ form a Steiner $(3,6,22)$ system. Pick 3 points $\{\alpha_1,\alpha_2,\alpha_3\}$ in the 22 points, a block $O_1$ of $\mathcal{S}(5,8,24)$ that contains $\{a,b,\alpha_1,\alpha_2\}$ , another block $O_2$ that contains $\{a,b,\alpha_3\}$ and intersect the first block at only 2 points $\{a,b\}$, and a third block $O_3$ that contains $\{a,\alpha_1,\alpha_2,\alpha_3\}$ and intersect both $B_1$ and $B_2$ at 4-point sets. It turns out that $M:=Stab(a)\cap Stab(b)\cong M_{22}$ and $M\cap Stab(\alpha_1), M\cap Stab(\alpha_2), M\cap Stab(\alpha_3), M\cap Stab(B_1),M\cap Stab(B_2),M\cap Stab(B_3)$ are maximal in $M$ and are in general position.\par

Sophie Le has checked computationally that $MaxDim(M_{22})=MaxDim(M_{23})=6$.\par
In a Steiner (5,8,24) system, a block is termed as an "octad". The stabilizer of an octad in $M_{24}$ is $C_2^4:A_8$, where $A_8$ permutes the 8 points in the octad and each involution in $C_2^4$ moves all 16 points outside the octad and stabilizes the points inside the octad \cite{CoSl88}. Thus we can use point stabilizers again to give a lower bound of $MaxDim$. Pick an octad in the Steiner system, then pick 6 points $a_1,...,a_6$ inside the block and one point $a_7$ outside the block. The stabilizers of them are then in general position. We shall prove this by showing that on one hand, the point-wise stabilizer of the 7 points is trivial, and on the other hand, point-wise stabilizers of any 6 points are either $C_3$ or $C_2^4$.\par
Since any five points determine a unique octad, a sixth point can either sit inside or outside the octad. Let the stabilizer of that octad be denoted $O\cong E:A$, where $E\cong E_{16}$ is elementary abelian of order 16 and $A\cong A_8$ is alternating on 8 points. Let $\pi_A:O\rightarrow A_8$ be the projection onto the permutation group of points in the octad, $\iota_E:E_{16}\rightarrow O$ be the injection from the normal subgroup, $\iota_A:A_8\rightarrow O$ be the splitting homomorphism, and the point-wise stabilizer of 6 points be denoted as $I$. If the sixth point sits inside the octad, then since any involution in the normal subgroup $E$ stabilizes the points in the octad, $\iota^{-1}_E(I)\cong E_{16}$. Since the nontrivial permutation on two points is odd, $\pi_A(I)$ must be trivial, thus $I\cong E_{16}$. If the sixth point sits outside the octad, then all involutions in $E$ cannot be contained in $I$, and there is a permutation $p$ in $A_8$ that shuffles the three points left over in the octad. Since $O$ is a split extension, $\iota_A(p)$ is an order 3 element in $I$. Since $\pi_A(I)\cong C_3$, $\iota^{-1}_E(I)$ is trivial, $I$ must be isomorphic to $C_3$. In either case, six-fold intersections of point stabilizers are non-trivial. However, if $R$ is the intersection of all 7 point stabilizers, then $\iota^{-1}_E(R)\cong \pi_A(R)\cong \mathbb{1}$, so $R$ must be trivial. Thus the point stabilizers are in general position and $MaxDim(M_{24})\geq 7$.\par
To obtain an upper bound, we checked computationally that all maximal subgroups of $M_{24}$ have $i\leq 6$ except the dyad groups $M_{22}:2$, the octad groups $C_2^4:A_8$ and the trio groups $C_2^6:(PSL_3(2)\times S_3)$, all of which have $i=7$. Thus $i(M_{24})\leq 8$ and if there is a length 8 list $l$ of maximal subgroups in general position, then only dyad, octad, and trio subgroups may occur in that list.\par
The dyad subgroup is the stabilizer of a two point set (i.e. a dyad). Since $M_{24}$ is 5-transitive, the intersection of two dyad groups falls into two conjugacy classes: if the dyads intersect, then the intersection is isomorphic to $PSL_3(4)$, and if the dyads are disjoint, then the intersection is isomorphic to $C_2^5:A_5:C_2:C_2$. Both of them have $i\leq 5$, so they cannot be in a length 8 list of Maximal subgroups in general position. That is, in a length 8 list $l$, there can be at most 1 dyad stabilizer. In the following paragraphs, the high transitivity of $M_{24}$ will be used extensively to show that beside the potential dyad stabilizer, we can use at most 2 octad stabilizers and 4 trio stabilizers, thus $l$ cannot exist.\par
\begin{lemma}
Had there been a length 8 list $l$ of maximal subgroups in general position, then octad stabilizers can appear at most twice in $l$.
\end{lemma}
\begin{proof}
Two distinct octads in the Steiner (5,8,24) system $\mathcal{S}(5,8,24)$ may be disjoint or intersect at a set of size 2 or 4. Any pair of octads that intersect at a 2-set are congruent to one another pair under the action of $M_{24}$. To see this, we can denote the 2-point intersections in each pair $\{a,b\}$ and $\{a',b'\}$ respectively. The octads involved in each pair shall be denoted $O_1,O_2$ and $O_1',O_2'$. Then using the 5-transitivity of $M_{24}$ we can find an element $\sigma_1$ that sends $a$ to $a'$ and $b$ to $b'$. The point-wise stabilizer of $a'$ and $b'$ is $M_{22}$, a group of automorphisms of the Steiner (3,6,22) system induced by fixing $a'$, $b'$ in $\mathcal{S}(5,8,24)$. As Brouwer mentioned, if we construct a $M_{22}$ graph by connecting the 6-point blocks a Steiner (3,6,22) system when they are disjoint, then $M_{22}$ would act edge-transitively on the graph \cite{Br89}. That is to say, since the 6-point blocks are constructed by removing $\{a',b'\}$ from the 77 octads in $\mathcal{S}(5,8,24)$ containing $\{a',b'\}$, there is an element $\sigma_2$ in the point-wise stabilizer of $a'$ and $b'$ that maps $\sigma_1(O_1)$ to $O_1'$ and  $\sigma_1(O_2)$ to $O_2'$. Thus $\sigma_1\sigma_2$ maps one pair of octad to the other pair, and it implies that $Stab(O_1)$ (resp.$Stab(O_2)$) can be mapped to $Stab(O_1')$ (resp.$Stab(O_2')$) under conjugation by $\sigma_1\sigma_2$. As a consequence, we only need to check one pair of octad stabilizers to know the intersection of any pair $(Stab(O_1),Stab(O_2))$ where $|O_1\cap O_2|$=2. It turns out that the intersection is isomorphic to $S_6$, which has $i=5$. Thus in $l$ two octad stabilizers must have their octads either disjoint or intersecting at a 4-point set.\par

Moreover, if an octad is split into a 4+2+2 partition by the blocks of a trio (i.e. a partition of the 24 points into 3 octads), then the intersection of the octad stabilizer with the trio stabilizer is $(C_2\times C_2\times A_4):2$, which has $i=4$, so they cannot coexist in $l$. As in the previous case, by using transitivity we only need to check one example in order to cover the entire case. Indeed, every octad-trio pair with this pattern of intersection are congruent to each other under the action of $M_{24}$, so the intersections of their stabilizers are all isomorphic. This claimed can be inferred from the 2+1+1 transitivity of the trio group \cite{CoSl88}. To prove the claim, denote the octad-trio pairs as $(O,T)$ (resp. $O',T'$) and the 3 octad blocks in $T$ as $O_1,O_2,O_3$ (resp. $O_1'$ through $O_3'$), where $|O_1\cap O|=4$. Let $\{a_1,a_2\}$ (resp. $a_1',a_2'$) be the two points in $O\cap O_2$ (resp. $O\cap O_2'$) and $a_3$ (resp. $a_3'$) be in $O\cap O_1$ (resp. $O\cap O_1'$), $a_4$ (resp. $a_4'$) be in $O\cap O_3$ (resp. $O\cap O_3'$). Since all trio groups are conjugate to each other, there is a $\sigma_1$ that maps trio $T$ to $T'$. Since a trio group acts as $S_3$ on its three octads, there is an element $\sigma_2$ that maps octad $\sigma_1(O_1)$ to $O_1'$ and $\sigma_1(O_2)$ to $O_2'$. Since the trio group is 2+1+1 transitive after stabilizing each octad in the trio individually, there is a $\sigma_3$ that takes $\sigma_2\sigma_1(\{a_i\})$ to $\{a_i'\}$. Now we have $\sigma_3\sigma_2\sigma_1(O\cap O_2)=O'\cap O_2'=\{a_1',a_2'\}$. If  $\sigma_3\sigma_2\sigma_1(O\cap O_3)=O'\cap O_3'$, then since $(O'\cap (O_2'\cup O_3'))\cup \{a_3'\}$ determines a unique octad, we have $\sigma_3\sigma_2\sigma_1(O)=O'$. If not, then $I:=\sigma_3\sigma_2\sigma_1(O\cap O_3)\cap(O'\cap O_3')$ must have only 1 element to ensure that $|O'\cap \sigma_3\sigma_2\sigma_1(O)|$ is less than 5. But in that case, the symmetric difference octad $O'\ominus O_1'$ would intersect $\sigma_3\sigma_2\sigma_1(O)$ in a 6 point set $(\sigma_3\sigma_2\sigma_1(O)\backslash \{a_3\})\cup \{a_1,a_2,a_4\}$, contradiction. Thus $\sigma_3\sigma_2\sigma_1(O)=O'$ and $\sigma_3\sigma_2\sigma_1$ maps the first octad-trio pair to the second.\par

If three octad stabilizers are in general position then their octads cannot be mutually disjoint, for by stabilizing two disjoint octads we automatically stabilize the third one. Thus there must be a pair of octad stabilizers intersecting at a tetrad (i.e. 4-point set). They will be denoted $O_1$ and $O_2$. In $\mathcal{S}(5,8,24)$, stabilizing a tetrad implies stabilizing a sextet (i.e. a partition of the 24 points into 6 tetrads). The tetrads in the sextet would be denoted as $Tet_1,...,Tet_6$, where $Tet_1:=O_1\cap O_2$, $Tet_2:=O_1-O_2$, and $Tet_3:=O_2-O_1$. Adding a third octad stabilizer into the list either preserves the tetrads in the sextet or refines them. However, the third octad $O_3$ must not intersect any tetrad in odd number of points. If $|O_3\cap Tet_i|=1$ for some $i$, then since whenever $i\neq j$, $Tet_i\cup Tet_j$ is octad and the intersection of any two distinct octads cannot have size other than 0, 2, or 4, $|O_3\cap Tet_j|=1$ or $3$ for all $i\neq j$. Because of the pigeon hole principal, we can then consider only the case where $|O_3\cap Tet_i|=3$ for some $i$ without loss of generality. If $|O_3\cap Tet_i|=3$ for some $i$, then $|O_3\cap Tet_j|=1$ for all $i\neq j$. Since we need $|O_3\cap O_1|=|O_3\cap O_2|=4$, $i$ must be 1. But then we have stabilized 2 points inside and 1 point outside $O_1$, namely we have stabilized $Tet_1-O_3$, $Tet_2\cap O_3$, and $Tet_3\cap O_3$. Thus the intersection of such three octad stabilizers must be isomorphic to a subgroup of $A_6$, and by the flatness of alternating groups it has $i\leq 4$ \cite{CaCa02}. If we want the list $l$ to have length 8, any three elements of it must have an intersection with $i\geq 5$, so $|O_3\cap Tet_i|=2$ or $4$ for all $i$.\par
If the third octad $O_3$ splits one of the tetrads into two dyads, then it must split $Tet_1,...,Tet_3$ as well as one more tetrad $Tet_4$ into dyads. This is because the 8 points in $O_3$ must be split into 4 dyads, so it must intersect with $O_1\cup O_2$. Since it must not intersect $O_1$ or $O_2$ in only 2 points, its intersection with $Tet_1,...,Tet_3$ must all have size 2. Moreover, the sextet induced by $|O_3\cap O_1|$ consists of tetrads that their intersection with $Tet_1,...,Tet_6$ could only have size 0 or 2. Let the tetrads in this sextet be denoted as $Tet_1',...,Tet_6'$, where $Tet_1':=O_1\cap O_3$, $Tet_2':=O_1-O_3$, and $Tet_3':=O_3-O_1$. The tetrad containing $Tet_3-O_3$ shall be denoted as $Tet_4'$. We know $Tet_4'=Tet_3-O_3\cup Tet_4-O_3$, because otherwise since there is an $i$ that $Tet_i'\subset Tet_5\cup Tet_6$, $(Tet_4'\cup Tet_i')\cap (Tet_5\cup Tet_6)$ would have cardinality 6, which is impossible. Therefore, $Tet_5\cup Tet_6=Tet_5'\cup Tet_6'$. Now we shall show that $|Stab(O_1)\cap Stab(O_2)\cap Stab(O_3)|\leq 16$ using those information. Pick a point $p$ in $Tet_5\cap Tet_5'$ and another point $q$ in $Tet_4\cap Tet_4'$ and denote $Stab(p)\cap Stab(q)\cap Stab(O_1)\cap Stab(O_2)\cap Stab(O_3)$ by $I$. Under the action of $Stab(O_1)\cap Stab(O_2)\cap Stab(O_3)$, $q$ can only be map to the other point in $Tet_4\cap Tet_4'$ and $p$ can only be mapped to points in the octad $Tet_5\cup Tet_6$, so $I$ must have index less than 16 in $Stab(O_1)\cap Stab(O_2)\cap Stab(O_3)$. But $I$ must be trivial, for by stabilizing $p$ we have stabilized all 12 dyads that are the intersections of $Tet_i\cap Tet_j'$, therefore every nontrivial element in $I$ must be involutions. Since we have stabilized an octad $Tet_5\cup Tet_6$, $I$ must act as an even permutation on the octad. If there is a nontrivial element $\sigma$ in $I$, then it must permute the two points in $Tet_5-Tet_5'$ or the two points in $Tet_5'-Tet_5$. In fact, they must permute both of them simultaneously to ensure that $\sigma$ is an even permutation on both $Tet_4\cup Tet_5$ and $Tet_4\cup Tet_6$. As a result $\sigma$ must stabilize the two points in $Tet_6\cap Tet_6'$ and $Tet_3\cap Tet_3'$ and move the points in $Tet_4\cap Tet_3'$. But an involution in $M_{24}$ moves either 16 or 24 points \cite{CoSl88}, so $\sigma$ must move every point in $O_1$, thus its action on $Tet_1\cup Tet_4$ is an odd permutation, contradiction. Since $i(I)$ is less than the number of prime factors in $I$, $i(I)\leq 4$, so we must avoid this pattern as well.\par

As a consequence, if three octad stabilizers are in general position, then their pairwise intersections must have cardinality 4, so $O_3$ is the union of two tetrads in $\{Tet_1,...,Tet_6\}$. We can add at most one more octad stabilizer into the list, for four octad stabilizers in general position would have stabilized all six tetrads individually. Meanwhile, we can add at most one trio stabilizer into the list. This is because on one hand, if the octads inside the trio are all unions of tetrads in $\{Tet_1,...,Tet_6\}$, then $Tet_5$, $Tet_6$ must be in different octads, so adding that trio stabilizer stabilizes all the 6 tetrads, and no more such trio stabilizers can be added. 
On the other hand, we can consider the case where one of the three octads in trio $\tilde{O}_1,\tilde{O}_2,\tilde{O}_3$, say $\tilde{O}_1$, intersect with a tetrad at a set of cardinality 2, or 3. Intersecting with a tetrad at one point implies intersecting with another one at a set of cardinality 3, so without loss of generality we can omit the discussion of that case. If $|\tilde{O}_1\cap Tet_i|=3$, then to ensure the intersection of $\tilde{O}_1$ with $O_1,O_2,O_3$ all have cardinality 0 or 4, we must have $i=1$ and $Tet_1\subset O_3$. But if we denote the octad containing $Tet_1-\tilde{O}_1$  as $\tilde{O}_2$, then $\tilde{O}_2\cap Tet_2$ and $\tilde{O}_2\cap Tet_3$ must all have cardinality 3 to ensure that $|\tilde{O}_2\cap O_1|=|\tilde{O}_2\cap O_2|=4$, so $|\tilde{O}_2\cap (Tet_2\cup Tet_3)|=6$, contradiction. If $|\tilde{O}_1\cap Tet_i|=2$, then again using the same pigeon hole argument we can infer that $|\tilde{O}_1\cap Tet_i|=|\tilde{O}_2\cap Tet_i|=2$ for $i=1,2,3,4$, where $Tet_4\subset O_3$. But then an arbitrary element $\sigma$ in the intersection of $Stab(O_1), Stab(O_2)$ and the trio stabilizer would not be able to move an element in $O_3$ outside, so the three octad stabilizers and the trio stabilizer are not in general position. Therefore, the only trio stabilizers that can be added in consist of octads that are unions of tetrads in $\{Tet_1,...,Tet_6\}$. Thus if there are three octad stabilizers in the list $l$, then we can add at most 1 dyad stabilizer, 1 octad stabilizer and 1 trio stabilizer into the list (Actually, at most 1 octad stabilizer or 1 trio stabilizer along with a potential dyad stabilizer), which means the list has size less than 8. It then follows immediately that we can use at most two octad stabilizers in $l$.
\end{proof}
\begin{lemma}
Had there been a length 8 list $l$ of maximal subgroups in general position, then trio stabilizers can appear at most 4 times in $l$.
\end{lemma}
\begin{proof}
Let us pick a trio $\{\tilde{O}_1,\tilde{O}_2,\tilde{O}_3\}$ in the Steiner system. A second trio $\{\tilde{Q}_1,\tilde{Q}_2,\tilde{Q}_3\}$ may intersect with the first trio in one of the following four ways:
\begin{enumerate}[label=(\roman*)]
\item All of $\tilde{Q}_1,\tilde{Q}_2,\tilde{Q}_3$ are split into 4+2+2 partitions.
\item $\tilde{Q}_1,\tilde{Q}_2$ are split into 4+2+2 partitions, and $\tilde{Q}_3$ is split into two tetrads.
\item All of $\tilde{Q}_1,\tilde{Q}_2,\tilde{Q}_3$ are split into tetrads.
\item $\tilde{Q}_1$ overlaps with $\tilde{O}_1$, and $\tilde{Q}_2,\tilde{Q}_3$ are split into tetrads.
\end{enumerate}
In this section we shall denote a pair of trios as $Trio_1=\{\tilde{O}_1,\tilde{O}_2,\tilde{O}_3\}, Trio_2=\{\tilde{Q}_1,\tilde{Q}_2,\tilde{Q}_3\}$ and in the case we are discussing more than one pairs, the second pair as $Trio_1'=\{\tilde{O}_1',\tilde{O}_2',\tilde{O}_3'\}, Trio_2'=\{\tilde{Q}_1',\tilde{Q}_2',\tilde{Q}_3'\}$. The intersection of the two trio stabilizers shall be denoted as $I(Trio_1,Trio_2)$.\par
Let the type $(i)$ pairs be considered the first. Assume without loss of generality that $|\tilde{O}_i\cap \tilde{Q}_i|=|\tilde{O}_i'\cap \tilde{Q}_i'|=4$ for $i=1,2,3$. Since $\tilde{O}_i\cap \tilde{Q}_i$ can only be mapped to $\tilde{O}_j\cap \tilde{Q}_j$ for $i,j=1,2,3$, $I(Trio_1,Trio_2)$ acts as a subgroup of $S_3$ on $\{\tilde{O}_i\cap \tilde{Q}_i|i=1,2,3\}$. $(\tilde{O}_i\cap \tilde{Q}_i)\cup (\tilde{O}_j\cap \tilde{Q}_j)$ cannot be octads in order to ensure that $(\tilde{Q}_i-\tilde{O}_i)\cup (\tilde{O}_j\cap \tilde{Q}_j)\cap \tilde{O}_j$ has order less than 5, so if we find the stabilizers $Stab(Sex_1)$ (resp.$Stab(Sex_2)$) of the sextet containing $\tilde{O}_1\cap \tilde{Q}_1$ (resp.$\tilde{O}_2\cap \tilde{Q}_2$), then $Stab(Sex_1)\cap Stab(Sex_2)\cap I(Trio_1,Trio_2)$, which is the stabilizer of each $\tilde{O}_i\cap \tilde{Q}_i$, must have index less than 6 in $I(Trio_1,Trio_2)$. Let the tetrads in $Sex_1$ be denoted as $Tet_1,...,Tet_6$. We know that $|Tet_i\cap (\tilde{O}_2\cap \tilde{Q}_2)|$ can only be 0, 1, or 2 to ensure $|((\tilde{Q}_1- \tilde{O}_1)\cup Tet_i)\cap \tilde{O}_2|$=2 or 4. But if there is no $i$ that $|Tet_i\cap (\tilde{O}_2\cap \tilde{Q}_2)|=1$, then there must be two tetrads so that $Tet_i\cup Tet_j=(\tilde{O}_2\cap \tilde{Q}_2)\cup (\tilde{O}_3\cap \tilde{Q}_3)$ in order to ensure that $|((\tilde{O}_1- \tilde{Q}_1)\cup Tet_i)\cap \tilde{Q}_2|=|((\tilde{O}_1- \tilde{Q}_1)\cup Tet_j)\cap \tilde{Q}_2|=4$. In that case, $(\tilde{O}_2\cap \tilde{Q}_2)\cup (\tilde{O}_3\cap \tilde{Q}_3)$ forms an octad, which is impossible. Thus both $\tilde{O}_2\cap \tilde{Q}_2$ and $\tilde{O}_3\cap \tilde{Q}_3$ are split into blocks of size 2+1+1 by tetrads in $Sex_1$. In particular, there is a tetrad $Tet_4$ that intersect both $\tilde{O}_2\cap \tilde{Q}_2$ and $\tilde{O}_3\cap \tilde{Q}_3$ at 2-sets and two tetrads $Tet_5,Tet_6$ that intersect both $\tilde{O}_2\cap \tilde{Q}_2$ and $\tilde{O}_3\cap \tilde{Q}_3$ at a single point. Pick one point $p\in Tet_4\cap \tilde{O}_2\cap \tilde{Q}_2$ and another point $q\in Tet_5\cap \tilde{O}_2\cap \tilde{Q}_2$, then $Stab(p)\cap Stab(q)\cap Stab(Sex_1)\cap Stab(Sex_2)\cap I(Trio_1,Trio_2)$ must be trivial, for by doing so we have stabilized an octad $\tilde{O}_2$, 6 points inside ($\{p\}, \{q\}, \tilde{O}_2\cap Tet_4-\{p\}, \tilde{O}_2\cap \tilde{Q}_2\cap Tet_6, \tilde{O}_2\cap Tet_5-\{q\}, \tilde{O}_2\cap Tet_6-\tilde{Q}_2$) and one point outside (say, $\tilde{O}_3\cap \tilde{Q}_2\cap Tet_5$) the octad. Since $p$ can only possibly be mapped to the other point in $\tilde{O}_2\cap Tet_4-\{p\}$ and $q$ can only be mapped to the point in $\tilde{O}_2\cap \tilde{Q}_2\cap Tet_6$ by elements in $Stab(Sex_1)\cap Stab(Sex_2)\cap I(Trio_1,Trio_2)$, the trivial group must have index less than or equal to 24 in $I(Trio_1,Trio_2)$. But then $I(Trio_1,Trio_2)$ cannot have more than 4 prime factors, so $i(I(Trio_1,Trio_2))\leq 4$, such a pair of trios stabilizers cannot appear in a list of length 8.\par
A pair of trio stabilizers of type $(ii)$ cannot appear either. In fact, if $\{Trio_1,Trio_2\}$ and $\{Trio_1',Trio_2'\}$
are such trio pairs, then we can show that $I(Trio_1,Trio_2)$ is conjugate to $I(Trio_1',Trio_2')$. Assume without loss of generality that $\tilde{O}_1,\tilde{Q}_3$ (resp.$\tilde{O}_1',\tilde{Q}_3$') are split into two tetrads. Denote$\tilde{O}_1\cap \tilde{Q}_1,\tilde{O}_1-\tilde{Q}_1,\tilde{Q}_1-\tilde{O}_1, \tilde{Q}_2-\tilde{O}_1$ as $Tet_1,...,Tet_4$. and their counterparts in the other pair of trios as $Tet_1',...,Tet_4'$. Those tetrads must be in the same sextet. Moreover, the two tetrads left over in $Tet_1,...,Tet_6$ must split both $\tilde{O}_2\cap \tilde{Q}_3$ and $\tilde{O}_3\cap \tilde{Q}_3$ in to dyads. Pick $p_1,...,p_4$ so that $p_1\in \tilde{Q}_1-\tilde{O}_1=Tet_3$, $p_2,p_3\in \tilde{Q}_3\cap\tilde{O}_2\cap Tet_5$ and $p_4\in \tilde{Q}_3\cap\tilde{O}_2\cap Tet_6$. Pick $p_1',...,p_4'$ in a similar way for the other pair of trios. Since all sextet stabilizers are conjugate to each other and a sextet stabilizer is transitive on the 6 tetrads, there is in element $\sigma\in M_{24}$ that maps $Tet_i$ to $Tet_i'$. Since as Conway mentioned the stabilizer of every tetrad in a sextet stabilizer is still 2+1+1+0+0+0 transitive \cite{CoSl88}, there must be a $\tau$ that fixes each tetrad and maps $p_i$ to $p_i'$. Now we have $\tau\sigma(\tilde{Q}_i)=\tilde{Q}_i'$ and $\tau\sigma(\tilde{O}_1)=\tilde{O}_1'$. If $\tau\sigma(\tilde{O}_2)\neq\tilde{O}_2'$ then $\tau\sigma(\tilde{O}_2)\cap\tilde{O}_2'=\{p_1,...,p_4\}$, therefore $|\tau\sigma(\tilde{O}_2)\cap\tilde{O}_3'\cap\tilde{Q}_3'|=1$. But if we consider the sextet induced by $\tilde{O}_3'\cap\tilde{Q}_3'$ then we know there is an octad consisting of 4 points in $\tilde{O}_1'$ and the 4 points in $\tilde{O}_3'\cap\tilde{Q}_3'$. That octad would then intersect $\tau\sigma(\tilde{O}_2)$ at a single point, which is not possible. Thus $\tau\sigma(\tilde{O}_2)=\tilde{O}_2'$ and $\tau\sigma$ maps the first pair of trios to the second pair. Checking an instance, we know that in this case $I(Trio_1,Trio_2)$ must be isomorphic to $D_8\times D_8$, which has $i=4$, and we can therefore eliminate all trio stabilizer pairs of this type.\par
For type $(iii)$ pairs, we know that $\{\tilde{O}_1,\tilde{O}_2,\tilde{O}_3,\tilde{Q}_1,\tilde{Q}_2,\tilde{Q}_3\}$ are all unions of tetrads in a particular sextet. A third trio may either consist of unions of tetrads in that particular sextet, or contain an octad intersecting one of the tetrads at a 2-point set. Otherwise, if there is an octad intersecting a tetrad in a 1-point set or 3 point set, then that octad must be partitioned into blocks of size 3+1+1+1+1+1, thus it must intersect one of the octads in $Trio_1$ or $Trio_2$ in a 2-point set, therefore we can no longer obtain a length 8 list in the end. But if the third trio contains an octad intersecting one of the tetrads at a 2-point set, say $Tet_1:=\tilde{O}_1\cap\tilde{Q}_1$, then it must intersect $Tet_2:=\tilde{O}_1-\tilde{Q}_1, Tet_3:=\tilde{Q}_1-\tilde{O}_1$ at 2-point sets. Let $\tilde{O}_2$ be the octad containing $Tet_3$ and $\tilde{Q}_2$ be the octad containing $Tet_2$. Then the octad in the third trio must intersect both  $Tet_4:=\tilde{Q}_2-Tet_2,Tet_5:=\tilde{O}_2-Tet_3$ at 2-point sets, which is impossible-we only have 8 points in an octad. Thus the third trio must consist of unions of tetrads in the sextet induced by $Tet_i$. That is to say, the intersection of all trio stabilizers in the current list must contain the stabilizer of each individual tetrads in $\{Tet_i\}$, which is proven by Conway to be isomorphic to $C_2^6:C_3$. If we take $\tilde{O}_i$ as vertices and  $\tilde{Q}_i$ as edges, then we can construct a triangle by attaching an edge to two vertices if their corresponding structures satisfies $\tilde{Q}_i\subset \tilde{O}_j\cup \tilde{O}_k$. To say that the action of $I(Trio_1,Trio_2)$ on the 6 tetrads preserves the octads is just saying that it acts as automorphisms of that triangle, therefore we can see that $I(Trio_1,Trio_2)$ acts as s subgroup of $D_6\cong S_3$ on those tetrads, and it implies that the index of $C_2^6:C_3$ in $I(Trio_1,Trio_2)$ has less than or equal to 2 prime factors. Thus if in our list of length 8 there are $Trio_1, Trio_2$ of type $(iii)$ then besides them we can add at most 2 more trio stabilizers into our list.\par
Now we shall consider the last case that each pair of trio stabilizers in our list is of type $(iv)$. Let $Trio_1,Trio_2$ be one of such a pair, $\tilde{O}_3=\tilde{Q}_3$, and the sextet induced by $\tilde{O}_1\cap \tilde{Q}_1$ be denoted as $\{Tet_i\}$. If we introduce a third trio $Trio_3=\{\tilde{R}_1,\tilde{R}_2,\tilde{R}_3\}$, then there are two cases: in the first one, all octads in $Trio_3$ are unions of tetrads in $\{Tet_i\}$; in the second case, there is an octad in $Trio_3$ that intersect one of the tetrads at a 2-point set. Using the same argument as before we can ignore the case where an octad in $Trio_3$ intersect one of the tetrads at a 1-point or 3-point set. To finish our proof we shall start with the first case. If $\tilde{O}_3=\tilde{Q}_3=\tilde{R}_3$ then $I(Trio_1,Trio_2)$ would be contained in the intersection of three trio stabilizers. To see this we can use $Tet_1:=\tilde{O}_1\cap\tilde{Q}_1,Tet_2:=\tilde{O}_1\cap\tilde{Q}_2,Tet_3:=\tilde{O}_2\cap\tilde{Q}_1,Tet_4:=\tilde{O}_2\cap\tilde{Q}_2$ as vertices and construct an edge-colored rectangle, where the edges correspond to the octads containing a pair of tetrads in $Tet_1,...,Tet_4$ and they are colored according to whether they belong to $Trio_1$ or $Trio_2$. Elements in $I(Trio_1,Trio_2)$ are exactly those elements in $M_23$ that act as automorphisms of the rectangle, and Elements in $I(Trio_1,Trio_2,Trio_3)$ are those elements preserving the diagonals as well. But an automorphism of that colored triangle would automatically preserve the diagonals, so$I(Trio_1,Trio_2)\subset I(Trio_1,Trio_2,Trio_3)$, the trio stabilizers are not in general position. Therefore the octad shared by $Trio_1,Trio_3$ must be either $\tilde{O}_1$ or $\tilde{O}_2$. On the other hand, the octad shared by $Trio_2,Trio_3$ must be either $\tilde{Q}_1$ or $\tilde{Q}_2$, so $Trio_3$ contains one of $\tilde{O}_1,\tilde{O}_2$ and one of $\tilde{Q}_1,\tilde{Q}_2$ at the same time, which is impossible. Thus we cannot add in a $Trio_3$ that has type $(iv)$ intersections with both $Trio_1$ and $Trio_2$ while keeping the stabilizers in general position without letting the octads in $Trio_3$ split the tetrads. But in the second case, if there is an octad $\tilde{R}_1$ intersecting one of the tetrads at a 2-point set, then that octad must split one of $Tet_1,..,Tet_4$ into dyads. Since that octad is not equal to any octad in $Trio_1,Trio_2$, all of its intersections with octads in $Trio_1,Trio_2$ need to have 4 points, so that octad must be split into 4 dyads by the 4 tetrads $Tet_1,..,Tet_4$. As a consequence $\tilde{R}_3=\tilde{O}_3=\tilde{Q}_3$. If we want to introduce a fourth trio $\{\tilde{P}_1,\tilde{P}_2,\tilde{P}_3\}$ that follows this pattern, then $\tilde{P}_1,\tilde{P}_2$ must split all the dyads outside $\tilde{P}_3=\tilde{O}_3=\tilde{Q}_3=\tilde{R}_3$, for otherwise by observing the action of three-fold intersections on dyads we know that the 4-fold intersection is the same as a 3-fold intersection. Using the same argument of inspecting the action of 4-fold intersection on the intersection of those dyads again, we can show that there is no hope to introduce a $Trio_5$ with its stabilizer in general position with the first 4 stabilizers in this case, therefore we can use at most 4 trio stabilizers in the list.
\end{proof}
Since in our list we can use at most 1 dyad stabilizer, 2 octad stabilizers and 4 trio stabilizers, the length of our list cannot be 8, so $MaxDim(M_{24})=7$.
\end{proof}

\section{The Higman-Sims Group and McLaughlin Group}
In this section, we utilize the fact that $MaxDim(M_{22})=6$ to give a lower bound for the $MaxDim$ of two sporadic simple groups, $HS$ and $McL$.
\begin{prop}
$MaxDim(HS)\geq 6$, and $MaxDim(McL)\geq 6$
\end{prop}
\begin{proof}
Both the Higman-Sims group and the McLaughlin group can be defined in terms of automorphisms of strongly regular graphs. As mentioned by Brouwer, the Higman-Sims graph has a 1+22+77 construction. Take the points and blocks in $\mathcal{S}(3,6,22)$ as vertices and add in an extra vertex $v$. Connect $v$ to all points of $\mathcal{S}(3,6,22)$, connect a point to a block if the point is in the block, connect two blocks if they are mutually disjoint, and keep points mutually unconnected, then we get a strongly regular graph of parameter $(100,22,0,6)$ which is the Higman-Sims graph \cite{Br89}. A stabilizer $M$ of vertices in $HS$ is isomorphic to$M_{22}$. Each point (resp. block) stabilizer in $M_{22}$ then correspond to vertex stabilizers from the 22-orbit (resp. 77-orbit) of $M$. It turns out that if we construct $\mathcal{S}(3,6,22)$ from $\mathcal{S}(5,8,24)$ by stabilizing 2 points, then all stabilizers in $M_{22}$ of those structures on $\mathcal{S}(3,6,22)$ induced by octads containing only one of the 2 points correspond to the intersections of $M$ with the maximal subgroups in one of the conjugacy classes of $PSU_3(5)$ in $HS$\cite{Br89}. Thus we can pick 6 maximal subgroups in general position in $M$ as described before and lift them to 5 vertex stabilizers and a $PSU_3(5)$ in $HS$, and the result list would be a length 6 list of maximal subgroups in general position for $HS$. \par
This is not the only way to find 6 maximal subgroups in general position in $HS$. In particular, if we pick 5 vertices in Higman-Sims graph that induces a pentagon and another vertex not adjacent to any of them, then their point stabilizers are in general position.\par
As mentioned by Brouwer again, the McLaughlin graph has a 22+77+176 construction.  Pick 2 points $\{a,b\}$ in Steiner $(5,8,24)$ system, there are 22 points $\{p\}$ leftover, 77 octads $\{O\}$ incident to both $\{a,b\}$ and 176 octads $\{Q\}$ incident to $a$ but not $b$. Take $\{p, O, Q\}$ as vertices, connect $p$ to $O$ if $p$ is not in $O$, connect $p$ to $Q$ if $p$ is in $Q$, connect $O$ to $O'$ if their intersection is $\{a,b\}$, connect $Q$ to $Q'$ if their intersection has 2 points, and connect $O$ to $Q$ if their intersection has size 4, then we get a strongly regular graph of parameter $(275,112,30,56)$ which is the McLaughlin graph \cite{Br89}. The stabilizer of a 22+77+176 partition is a maximal subgroup $M$ of $McL$ isomorphic to $M_{22}$. By the construction of the graph, we can pick 3 vertices from the 22-orbit, 2 vertices from the 77 orbit, and 1 point from the 176 orbit such that the relation between their corresponding structures are exactly the same as those corresponding to the length 6 general position sequence in $M_22$ described above. Their stabilizers in $M$ are therefore in general position. Since each vertex stabilizer in $M$ can be lifted to a maximal subgroup in $McL$, we have found a list of maximal subgroups in general position, so $MaxDim(McL)\geq 6$.
\begin{remark}
The 176 octads $\{Q\}$ incident to $a$ but not $b$ (resp. to $b$ but not $a$) are the points (resp. quadrics) in Higman geometry \cite{Br89}.
\end{remark}
\end{proof}

\section{Remarks}
\begin{remark}
By looking at intersections of point stabilizers, we can show that $MaxDim(Sp_n(q))\geq n$ and $MaxDim(U_n(q))\geq n$. As Brouwer pointed out, the symplectic groups and unitary groups act on strongly regular graphs. However, this lower bound is usually not idealistic.
\end{remark}
Every group with $MaxDim>m$ observed so far has $i>MaxDim$ as well. Is it necessarily true for all groups with $MaxDim>m$? In particular, does it contradict the following conjecture of Shareshian? 
\begin{conj}
(Shareshian) Any open interval of a subgroup lattice of some finite group is has the same homotopy as a wedge of spheres. 
\end{conj}

\bibliographystyle{amsplain}
\bibliography{citations}

\end{document}